\documentclass[10pt]{article}
\usepackage{amsmath,amssymb,amsbsy,amsfonts,latexsym,amsthm,amscd,amsxtra,
amssymb}
\usepackage[all]{xy}
\usepackage{pdftricks}
\usepackage{epstopdf}

\newtheorem{Teo}{Theorem}%[section]
\newtheorem{Def}{Definition}
\newtheorem{Prop}{Proposition}
\newtheorem{Cor}{Corollary}
\newtheorem{Lem}{Lemma}
\newtheorem{Rema}{Remark}
\newenvironment{Rem}{\begin{Rema} \begin{upshape}} {\end{upshape}\end{Rema}}
\newtheorem*{Pf}{Proof}
\newenvironment{Proof}{\begin{Pf} \begin{upshape}} {\end{upshape} \qed\end{Pf}}

\newcommand\beqa[1]{ \begin{eqnarray} \label{#1}}
\newcommand{\eeqa}{ \end{eqnarray} }
\newcommand{\beqano}{ \begin{eqnarray*} }
\newcommand{\eeqano}{ \end{eqnarray*} }

\newcommand{\T}{ {\mathbb T}   }
\newcommand{\N}{ {\mathbb N}   }
\newcommand{\R}{ {\mathbb R}   }
\newcommand{\Z}{ {\mathbb Z}   }

\newcommand{\K}{ {\mathbb K}   }
\renewcommand \a {\alpha}
\newcommand \e {\varepsilon }

\renewcommand \b  {\beta}

\renewcommand \d {\delta}

\newcommand \m {\mu}

\newcommand \g {\gamma}

\renewcommand \l {\lambda}
\renewcommand \L {\Lambda}

\newcommand \cA {{\mathcal A}}

\newcommand \cF {{\mathcal F}}
\newcommand \cH {{\mathcal H}}
\newcommand \cL {{\mathcal L}}
\newcommand \cM {{\mathcal M}}

\newcommand \cS {{\mathcal S}}

\newcommand \calM {\mathfrak{M}}

\def\ie{\hbox{\it i.e.\ }}

\newcommand \dpr {\partial}

\newcommand \rH {{\rm H}}
\newcommand \rT {{\rm T}}

\begin{document}

\title{Differentiability of Mather's average action and integrability on closed surfaces}
\author{\sc{Daniel Massart and Alfonso Sorrentino}}
\date{}
\maketitle

%%%%%%%%%%%%%%%%%%%%%%%%%%%%%%%%%%%%%%%%%%%%%%%%
%%%%%%%%%%%%%%%%%% ABSTRACT %%%%%%%%%%%%%%%%%%%%
\begin{center}
\parbox[c]{4.5 in}
{\small In this article we study the differentiability of Mather's average action, or $\beta$-function, on closed surfaces, with particular attention to its relation to the integrability of the system.}\\
\end{center}

%\vspace{5pt}

%%%%%%%%%%%%%%%%%%%%%%%%%%%%%%%%%%%%%%

\section{Introduction}\label{section1}

In the study of Tonelli Lagrangian and Hamiltonian systems,  a central role in understanding  the dynamical and topological properties of the  action-minimizing sets (also called {\it Mather} and {\it Aubry sets}), is played by the so-called {\it Mather's average action} (sometimes referred to as $\beta$-{\it function}  or {\it effective Lagrangian}), with particular attention to its differentiability and non-differentiability properties. %and its convex conjugate $\alpha$ function.
Roughly speaking,  this is a convex superlinear function on the first homology group of the base manifold, which represents the minimal action of invariant probability measures within a prescribed  {\it homology class}, or {\it rotation vector} (see (\ref{defbeta}) for a more precise definition). 
Understanding whether or not this function is differentiable, or even smoother, %what happens at non-differentiability points, or 
and what are the implications of its regularity to the dynamics  of the system is a formidable problem, which is still far from being completely understood.
Examples of Lagrangians  admitting a smooth $\beta$-function are easy to construct; trivially, if the base manifold $M$ is such that $\dim \rH_1(M;\R)=0$ then $\beta$ is a function defined on a single point and it is therefore smooth. Furthermore, if $\dim \rH_1(M;\R)=1$ then a result by Carneiro \cite{Carneiro}  allows one to conclude that $\beta$ is differentiable everywhere, except possibly at the origin. 
As soon as  $\dim \rH_1(M;\R)\geq 2$ the situation becomes definitely less clear and the smoothness of $\beta$ becomes a more ``untypical'' phenomenon. Nevertheless, it is still possible to find some interesting examples in which it is smooth.
 For instance, let 
$H:\rT^*\T^n \longrightarrow \R$ be a completely integrable Tonelli Hamiltonian system, given by $H(x,p)=h(p)$, and consider the associated Lagrangian $L(x,v)=\ell(v)$ on $\rT \T^n$.  It is easy to check that in this case, up to identifying $\rH_1(\T^n;\R)$ with $\R^n$, one has $\beta(h)=\ell(h)$ and therefore $\b$ is as smooth as the Lagrangian and the Hamiltonian are.  One can weaken the assumption on the complete integrability of the system and consider $C^0$-{\it integrable systems}, \ie Hamiltonian systems that admit a foliation of the phase space by disjoint invariant continuous Lagrangian graphs, one for each possible cohomology class (see Definition \ref{C0integrability}). It is then possible to prove that also in this case the associated $\beta$ function is $C^1$ (see Lemma \ref{lemma2}). 
These observations raise the following question: with the exception of the mentioned trivial cases, does the regularity of $\b$ imply the integrability of the system?

In this article we address the above problem in the case of Tonelli Lagrangians on closed surfaces, not necessarily orientable. The main results can be summarized as follows.\\

\noindent{\bf Main Results}. {\it
Let $M$ be a closed surface and $L:\rT M \longrightarrow \R$ a Tonelli Lagrangian. 
\begin{itemize}
\item[{\rm (i)}]   If $M$ is  not the sphere,  the projective plane, the Klein bottle or the torus, then $\beta$ cannot be $C^1$ on the whole of $\rH_1(M;\R)${\rm [Corollaries  \ref{oriented, not differentiable} and \ref{maintheo2, nonor}].}
\item[{\rm (ii)}] If $M$ is  not the torus then the Lagrangian cannot be $C^0$-integrable {\rm [Proposition \ref{torus_only_integrable}].}
\item[{\rm (iii)}]	If $M$ is the torus, then $\beta$ is $C^1$ if and only if the system is $C^0$-integrable {\rm [Theorem \ref{maintheo1}].}
\end{itemize}
}
\noindent In particular, in Corollary \ref{corollarymaintheo} we deduce several information on the dynamics of $C^0$-integrable systems. Moreover in section \ref{mechanicalcase} we discuss the case of mechanical Lagrangians on the two-torus and  show that  in this case: $\beta$ {\it is $C^1$ if and only if the  system is the geodesic flow associated to a flat metric on $\T^2$} [Proposition \ref{prop_mechanical_2dim}].\\

The paper is organized as follows. In Section \ref{section2} we provide a brief overview of Mather's theory for Tonelli Lagrangian systems and introduce the minimal average action, or $\beta$ function. In Section \ref{section3} we discuss several properties of the $\beta$ function associated to Lagrangians on closed surfaces, with particular attention to the implications of its differentiability, or the lack thereof, to the structure of the Mather and Aubry sets. We consider both orientable and non-orientable surfaces, pointing out the differences between the two cases. Finally in Section \ref{section4} we introduce  the concept of $C^0$-integrability and relate it to the regularity of $\beta$.

\vspace{15 pt}

\noindent{\bf Acknowledgements.} 
The authors are particularly grateful to Patrick Bernard, Victor Bangert and Albert Fathi for several useful discussions.
This work has been partially funded by the ANR project  ``{\it Hamilton-Jacobi et th\'eorie KAM faible}''  
(ANR-07-BLAN-3-187245). A.S. also wishes to acknowledge the support of {\it Fondation Sciences Math\'ematiques de Paris} and {\it Herchel-Smith foundation}.

%%%%%%%%%%%%%%%%%%%%%%%%%%%%%%%%%%%%%

\section{Mather's average action and Aubry-Mather theory}\label{section2}

Let us start by recalling some basic facts about Mather's theory for Tonelli Lagrangians. We refer the reader to \cite{Mather91, Mather93} for the proofs of these results (see also \cite{SorLecturenotes}).
Let $M$ be a compact and connected smooth $n$-dimensional manifold without boundary. Denote by ${\rm T}M$ its tangent bundle and ${\rm T}^*M$ the cotangent one and denote points in ${\rT M}$ and ${\rT^*M}$
respectively by $(x,v)$ and $(x,p)$. We will also assume that the cotangent bundle $\rT^*M$ is equipped with the canonical symplectic structure, which we will denote $\omega$.
A {\it Tonelli Lagrangian} is a $C^2$ function $L: \rT M \rightarrow \R$, which is strictly convex and uniformly superlinear in the fibers; in particular this Lagrangian defines a flow on ${\rT M}$, known as {\it Euler-Lagrange flow} $\Phi^L_t$, whose integral curves are solutions of $\frac{d}{dt}\frac{\dpr L}{\dpr v}(x,v)= \frac{\dpr L}{\dpr x}(x,v)$. 
Let $\calM(L)$ be the space of compactly supported probability measures on 
${\rm T}M$  invariant  under the Euler-Lagrange flow of $L$. To every 
$\m \in \calM(L)$, we may associate its {\it average action}
$$ A_L(\m) = \int_{{\rm T}M} L\,d\m\,. $$
It is quite easy to check that since 
$\m$ is invariant by 
the Euler-Lagrange flow, then for each $f\in C^1(M)$ we have 
$\int{{df}(x)\cdot v }\, d\m =0$ (see for instance \cite[Lemma on page 176]{Mather91}). Therefore we can define a linear functional
\beqano
{\rm H}^1(M;\R) &\longrightarrow& \R \\
c &\longmapsto& \int_{{\rm T}M} {\eta_c}(x)\cdot v\, d\m\,,
\eeqano
where $\eta_c$ is any closed $1$-form on $M$, whose cohomology class is $c$. By 
duality, there 
exists $\rho (\m) \in {\rm H}_1(M;\R)$ such that
$$
\int_{{\rm T}M} {\eta_c}(x)\cdot v\, d\m = \langle c,\rho(\m) \rangle
\qquad \forall\,c\in {\rm H}^1(M;\R)$$ 
(the bracket on the right--hand side denotes the canonical pairing between 
cohomology and 
homology). We call $\rho(\m)$ the {\it rotation vector} of $\m$. 
It is possible to show \cite[Proposition on page 178]{Mather91} that the map $\rho: \calM(L) \longrightarrow {\rm H}_1(M;\R)$ is surjective and hence there exist invariant probability measures for each rotation vector.
Let us consider the minimal value of the average action $A_L$ over the set of 
probability measures with rotation vector $h$ (this minimum exists because of the lower semicontinuity of the action functional \cite[Lemma on page 176]{Mather91})):
\begin{eqnarray}
\b: {\rm H}_1(M;\R) &\longrightarrow& \R \nonumber\\
h &\longmapsto& \min_{\m\in\calM(L):\,\rho(\m)=h} A_L(\m)\,. \label{defbeta}
\end{eqnarray}
This function $\beta$ is what is generally known as {\it Mather's 
$\beta$-function}.
A measure $\m \in \calM(L)$ realizing such a minimum amongst all invariant probability measures with the same rotation vector, \ie $A_L(\m) = 
\b(\rho(\m))$, is called an {\it action minimizing} {\it measure with rotation vector} $\rho(\mu)$ or, equivalently, a $(L,\rho(\mu))$-{\it minimizing measure}. 
The $\beta$-function is convex, and therefore one can consider its {\it conjugate} 
function (given by Fenchel's duality)
$ \a: {\rm H}^1(M;\R) \longrightarrow \R $ defined by
\beqano
\a(c) &:=& \max_{h\in {\rm H}_1(M;\R)} \left(\langle c,h \rangle - \b(h) 
\right)
= - \min_{h\in {\rm H}_1(M;\R)} \left(\b(h) - \langle c,h \rangle\right) =\\
&=& - \min_{\m \in \calM(L)} \left( A_L(\m) - \langle c, \rho(\m)\rangle\right)
= - \min_{\m \in \calM(L)} A_{L-\eta_c}(\m),\\
\eeqano
where $\eta_c$ is any smooth closed $1$-form on $M$ with cohomology class $c$. Observe that  the modified Lagrangian $L-\eta_c$  ($\eta_c$ can be also seen as  a function on $\rT M$: $\eta_c(x,v):=\eta_c(x)\cdot v$)  is still of Tonelli type and has the same Euler-Lagrange flow as $L$ (it follows easily from the closedness of $\eta_c$). Nevertheless, invariant probability measures that minimize the action functional  $A_{L-\eta_c}$ change accordingly to the chosen cohomology class.
A measure  $\m\in\calM(L)$ realizing the minimum of  $A_{L-\eta_c}$ (withouth any constraint on the rotation vector) is called {\it $c$-action minimizing} or, equivalently, $(L-c)$-{\it minimizing measure}.
This leads to the definition of a first important family of invariant subsets of $\rT M$:

\begin{itemize}
\item for a homology class $h\in \rH_1(M;\R)$, we define  the {\it Mather set of 
rotation vector} $h$ as:
$$ \widetilde{\cM}^h := \bigcup \left\{ {\rm supp}\,\m: \; \m\; \mbox{is action minimizing with rotation vector }\; h\right\};
$$

\item for a cohomology class $c \in {\rm H}^1(M;\R)$, we define the 
{\it Mather set of 
cohomology class} $c$ as:
$$ \widetilde{\cM}_c := \bigcup \left \{ {\rm supp}\,\m: \; \m\; {\rm is }\;c-\mbox{action minimizing}\right\}.
$$
\end{itemize}

The relation between these sets is described in the following lemma. To state it, let us recall that, like any convex 
function on a finite-dimensional space, Mather's   $\b$ function 
admits a subderivative at any point $h\in \rH_1(M;\R)$, \ie we can find $c\in 
\rH^1(M;\R)$ such that
$$\forall h'\in \rH_1(M;\R), \quad \b(h')-\b(h)\geq \langle c,h'-h\rangle.$$
As it is usually done, we will denote by $\partial \b(h)$ the set of $c\in 
\rH^1(M;\R)$ which are subderivatives of $\b$ at $h$, \ie the set of $c$'s which satisfy the  
above inequality (this set is also called the {\it Legendre transform of} $h$). 
By Fenchel's duality,  we have
\beqa{fenchelequality}
c\in \partial \b(h) \Longleftrightarrow \langle c,h\rangle=\a(c)+\b(h).
\eeqa

\begin{Lem}\label{lemma1}
Let $h,c$ be respectively an arbitrary homology class in $\rH_1(M;\R)$ and an arbitrary cohomology class $ \rH^1(M;\R)$. We have
$$
{\rm (1)}\; \widetilde{\cM}^{h} \subseteq \widetilde{\cM}_c \quad 
\Longleftrightarrow \quad {\rm (2)}\;c \in \partial \b(h)\,.
$$
\end{Lem}

\begin{Proof}
 Let us prove that
 ${\rm  (1)} \Longrightarrow {\rm (2)}$.
If $\widetilde{\cM}^{h}  \subseteq  \widetilde{\cM}_c $, then there exists a 
$c$-action minimizing invariant measure $\m$ with rotation vector $h$. 
Let $\eta_c$ be a closed $1$-form with $[\eta_c]=c$. 
From the definitions of $\a$, $\b$ and rotation vector:
\beqano
-\a(c) = \int_{{\rm T}M} (L-{\eta_c})\,d\m = \int_{{\rm T}M} L\,d\m - 
\langle c,h\rangle = \b(h) - \langle c,h\rangle\,;
\eeqano
since $\b$ and $\a$ are convex conjugated,  then $c$ is a subderivative of 
$\b$ at $h$.\\
Finally, in order to show $(2) \Longrightarrow (1)$, let us prove that any 
action minimizing measure with rotation vector $h$ is $c$-action minimizing. 
In fact, if $c \in \partial \b(h)$ then $\b(h)= \langle c,h\rangle - 
\a(c)$; therefore for any $\m \in \calM(L)$ with $\rho(\mu)={h}$ and $\eta_c$ as above:
$$
-\a(c) = \b(h) - \langle c,h\rangle =  \int_{{\rm T}M}(L-{\eta_c})\, 
d\m.  
$$
This proves that $\m$ is $c$-action minimizing and it concludes the proof.
\end{Proof}

\begin{Rem}\label{mathersetsdifferentcohomology}
One can also find a relation between the Mather sets corresponding to different cohomology classes, in terms of the strict  convexity of the $\a$-function, or better the lack thereof (see also  Section \ref{section2}).
The regions where $\a$ is not strictly convex are called {\it flats}: for instance the Legendre transform of $\a$ at $c$, denoted $\dpr \a(c)$, which is the set of  homology classes for which equality (\ref{fenchelequality}) holds, is an example of flat. It is possible to check that if two cohomology classes are in the relative interior of the same flat $F$ of $\a$, then their Mather sets coincide (see \cite{Mather91, Massart03}). We denote by $\widetilde{\cM}(F)$ the common Mather set to all the cohomologies in the relative  interior of $F$. 
\end{Rem}

In addition to the Mather sets, one can also construct another interesting family of compact invariant sets. We define the {\it  Aubry set $\widetilde{\cA}_c$ of cohomology class $c$} by looking at a special kind of global minimizers: we say that $(x,v) \in \rT M$ lies in $\widetilde{\cA}_c$ if there exists a sequence of $C^1$ curves $\gamma_n : \left[0,T_n\right] \longrightarrow M$, such that \begin{itemize}
  \item $\gamma_n (0)= \gamma_n (T_n) = x$ for all $n$
  \item $\dot{\g}_n(0)$ and $\dot{\g}_n(T_n)$ tend to $v$ as $n$ tends to infinity
       \item $T_n$ goes to infinity with $n$
  \item $\int_{0}^{T_n}  (L-\eta_c + \a(c))(\g_n(t),\dot{\g}_n(t))\,dt $ tends to $0$ as $n$ tends to infinity.
  \end{itemize}
One can prove \cite{Mather93, Fathibook} that if $(x,v)\in \widetilde{\cA}_c$, then the curve $\g(t):=\pi \Phi^L_t$, where $\pi: \rT M \longrightarrow M$ denotes the canonical projection, is a $c$-global minimizer. However, not all $c$-global minimizers can be obtained in this way.

\noindent These action-minimizing sets that we have defined, are such that $\widetilde{\cM}_c \subseteq \widetilde{\cA}_c$ for all $c\in\rH^1(M;\R)$ and moreover one of their most important features is that they are graphs over $M$ ({\it Mather's graph theorem} \cite{Mather91, Mather93}), \ie 
the projection along the fibers $\pi|{\widetilde{\cA}_c(L)}$ is injective and its inverse $\big(\pi|\widetilde{\cA}_c(L)\big)^{-1}\!\!\!: \pi\big(\widetilde{\cA}_c(L)\big) \longrightarrow \widetilde{\cA}_c(L)$ is Lipschitz.
Hereafter we will denote by $\cM^h$, $\cM_c$ and $\cA_c$ the corresponding projected sets.\\

Another interesting characterization of the Aubry set is provided by {\it weak KAM theory} \cite{Fathibook} in terms of {\it critical subsolutions} of Hamilton-Jacobi equation or, in a more geometric way, in terms of some special {\it Lipschitz Lagrangian graphs}. Let us consider the Hamiltonian system associated to our Tonelli Lagrangian. In fact, if one considers the {Legendre transform} associated to $L$, \ie the diffeomorphism
$\cL_L: \rT M \longrightarrow \rT^*M$ defined by
$\cL_L(x,v)=(x, \frac{\dpr L}{\dpr v}(x,v))$, then it is possible to introduce a {Hamiltonian system} 
$H: \rT^* M \rightarrow \R$, where $H(x,p)=\sup_{v\in \rT_xM}(\langle p,v \rangle - L(x,v))$. It is easy to check that $H$ is also $C^2$, strictly convex and uniformly superlinear in each fiber: $H$ is also said to be {\it Tonelli} (or sometimes {\it optical}). Recall that the flow on $\rT^* M$ associated to this Hamiltonian, known as the {\it Hamiltonian flow} of $H$, is conjugated - via the Legendre transform - to the Euler-Lagrange flow of $L$.
Therefore one can define the analogue of the Mather and Aubry sets in the cotangent space, simply considering $\cM^*_c:=\cL_L\big(\widetilde{\cM}_c\big)$ and
$\cA^*_c:=\cL_L\big(\widetilde{\cA}_c\big)$. These sets still satisfy all the properties mentioned above, including the graph theorem.
Moreover, they will be contained in the energy level $\{H(x,p)=\a(c)\}$, where $\a: \rH^1(M;\R)\longrightarrow \R$ is exactly {\it Mather's $\a$-function} introduced before. \\
Using Fathi's weak KAM theory \cite{Fathibook} it is possible to obtain a nice characterization of the Aubry set in terms of {critical subsolutions} of Hamilton-Jacobi equation. As above, let $\eta_c$ be a closed $1$-form with cohomology class $c$; we will say that $u \in C^{1,1}(M)$ is an $\eta_c$-critical subsolution if it satisfies $H(x,\eta_c + d_xu)\leq \a(c)$ for all $x\in M$. The existence of such functions has been showed by Bernard \cite{bernardc11}.
If one denotes  by $\cS_{\eta_c}$ the set of $\eta_c$ critical subsolutions, then: 
\beqa{Aubry} {\cA}_c^*=\bigcap_{u\in \cS_{\eta_c}} \left\{(x,\eta_c(x) + d_xu):\; x\in M\right\}.\eeqa
This set does not depend on the particular choice of $\eta_c$, but only on its cohomology class.
Observe that since in $\rT^*M$, with the standard symplectic form, there is a $1$-$1$ correspondence between Lagrangian graphs and closed $1$-forms, then %(see for instance \cite{AnnaCannas}). 
 we can interpret the graphs of these critical subsolutions as { Lipschitz} Lagrangian graphs in $\rT^*M$ and the Aubry set may be defined as their  {\it non-removable} intersection (see also \cite{PPS}).

%%%%%%%%%%%%%%%%%%%%%%%%%%%%%%%%%%
\section{Differentiability of $\beta$ on closed surfaces}\label{section3}
%%%%%%%%%%%%%%%%%%%
In this section we have gathered all the material we need on the differentiability of $\beta$. Let $M$ be a closed surface, not necessarily orientable,  $L: \rT M \longrightarrow \R$ a Tonelli Lagrangian on $M$ and $H: \rT^* M \longrightarrow \R$ the associated Hamiltonian. 

Let us recall that a homology class $h$ is said to be $k$-irrational, if $k$ is the dimension of the smallest subspace of $\rH_1(M;\R)$ generated by integer classes and containing $h$. In particular, $1$-{\it irrational} means ``on a line with rational slope'', while {\it completely irrational} stands for $\dim \rH_1(M;\R)$-irrational. Moreover, a homology $h$ is said to be {\it singular} if its Legendre transform $\dpr \beta(h)$ is a {\it singular flat}, \ie its Mather set $\widetilde{\cM}(\dpr \beta(h))$ - see Remark \ref{mathersetsdifferentcohomology} -  contains fixed points. Observe that the set of singular classes, unless it is empty,  contains the zero class and is compact.

For $h \in \rH_1(M;\R)\setminus \left\{0\right\}$, we define the {\it maximal radial flat} $R_h$ of $\beta$ containing $h$ as the largest subset of the half-line $\left\{th \; : t \in \left[0,+\infty \right) \right\}$ containing $h$ (not necessarily in its relative interior) in restriction to which $\beta$ is affine. The possibility of radial flats is the most obvious difference between the $\beta$ functions of Riemannian metrics (see \cite{Massart97b, nonor}) and those of general Lagrangians. An instance of radial flat is found in \cite{Carneiro-Lopes}. 
%We define the Mather set of $R_h$  as  the closure in $\rT M$ of the union of the supports of all $(L,th)$-minimizing measures, for $th \in R_h$.

%\begin{Rem} We will sometimes denote the Mather set of $R_h$ by $\widetilde{\cM}(R_h)$. We warn the reader against confusing this set with the Mather set associated to a flat $F$ of $\alpha$, 
%$\widetilde{\cM}(F)$, introduced in Remark \ref{mathersetsdifferentcohomology}.
%\end{Rem}

Let $h$ be a homology class. Assume $h$ is  1-irrational.
Then for all $t$ such that $th \in R_h$, $th$ is also 1-irrational. So every $(L,th)$-minimizing measure is supported on periodic orbits (see \cite[Proposition 5.6]{nonor}). Furthermore $R_h$ is  a flat of $\beta$, so there exists some cohomology class $c$ such that 
for every $t$ such that $th \in R_h$, for every $(L,th)$-minimizing measure $\mu$, the support of $\mu$ is contained in the Mather set $\widetilde{\cM}_c$. Then Mather's Graph Theorem  says that the projections of the periodic orbits  that support the $(L,th)$-minimizing measures, for every $t$ such that $th \in R_h$,  are pairwise disjoint. We usually denote  these periodic orbits by $(\gamma_i, \dot\gamma_i)_{i \in I}$, where $I$ is  some (possibly infinite) set.  Recall that  \cite[Theorem 1.3 and Proposition 2.4]{Massart09} combine to say the following.

\begin{Teo}\label{AvsM, principal}
Assume
\begin{itemize}
	\item $M$ is a closed surface
	\item $L$ is an autonomous  Tonelli Lagrangian on $M$
	\item $h$ is a 1-irrational, non-singular homology class
	\item $c$ lies in the relative interior of $\partial \beta (h)$
	\item $(\gamma_i, \dot\gamma_i)_{i \in I}$ are the periodic orbits which comprise the supports of all $(L,th)$-minimizing measures, for all $th$ in $R_h$.  
\end{itemize}
Then, the Mather set $\widetilde{\mathcal{M}}_c$  is precisely the union of the $(\gamma_i, \dot\gamma_i)$, and $\widetilde{\mathcal{A}}_c=\widetilde{\mathcal{M}}_c$.
\end{Teo}

From this we deduce this lemma.

\begin{Lem}\label{h_i conv partial alpha}
Assume
\begin{itemize}
	\item $M$ is a closed surface
	\item $L$ is an autonomous  Tonelli Lagrangian on $M$
	\item $h_0$ is a 1-irrational, non-singular homology class
	\item $c$ lies in the relative interior of $\partial \beta (h_0)$
	\item $(\gamma_i, \dot\gamma_i)_{i \in I}$ are the periodic orbits which comprise the supports of all $(L,th)$-minimizing measures, for all $th$ in $R_h$
	\item $\mu_i$ is the probability measure equidistributed on $(\gamma_i, \dot\gamma_i)$, for any $i \in I$
	\item $h_i$ is the homology class of $\mu_i$,  for any $i \in I$. 
\end{itemize}
Then, there exists a finite subset $J$ of $I$ such that $\partial \alpha (c)$ is the convex hull of the $h_i$, $i \in J$.
\end{Lem}
\begin{Proof}
For any $i \in I$, $\mu_i$ is $(L-c)$-minimizing, so $h_i \in \partial \alpha (c)$. By the convexity of $\partial \alpha (c)$, it follows that the convex hull of 
the $h_i$, $i \in I$, is contained in $\partial \alpha (c)$. 

Conversely, take $h \in \partial \alpha (c)$ and an $(L,h)$-minimizing measure $\mu$. Then $\mu$ is $(L-c)$-minimizing because $h \in \partial \alpha (c)$, therefore the support of $\mu$ is contained in the Mather set $\widetilde{\mathcal{M}}_c$, that is, $\mu$ is supported by some or all of the $(\gamma_i, \dot\gamma_i)$. Hence the homology class of $\mu$, which is $h$, is a convex combination of the $h_i$.  So $\partial \alpha (c)$ is contained in the convex hull of  the $h_i$, $i \in I$. 

We still have to find the finite subset $J$. Recall that if $\alpha_j$, $j \in J$, is a collection of pairwise disjoint closed curves on a  compact surface $M$,
then the collection  of homology classes $\left[\alpha_j\right]$, $j \in J$, is finite. Moreover its cardinality is at most  $3/2 (\dim \rH_1 (M;\R)-1)$ if $M$ is orientable. See, for instance, \cite{nonor} for an upper bound on non-orientable surfaces.  

Since the closed curves $\gamma_i$ are pairwise disjoint, their homology classes form a finite subset $h_1, \ldots h_n$ of $\rH_1(M;\R)$. Now, for every $i \in I$, the homology class $h_i$ of $\mu_i$ is $T_i^{-1} \left[ \gamma_i \right]$, where $T_i$ is the period of the periodic orbit $(\gamma_i, \dot\gamma_i)$. So all the homology classes $h_i$, $i \in I$, lie  in the finite union of straight lines $\R h_1 \cup \ldots \cup \R h_n$. Thus  the convex hull of  the $h_i$, $i \in I$, has at most $2n$ extremal points, two in each straight line
$\R h_1 , \ldots , \R h_n$. Therefore $\partial \alpha (c)$, which is the convex hull of  the $h_i$, $i \in I$, is a finite polytope.
\end{Proof}

Now we can state the main results of this section.  
First let us get the Klein bottle case out of the way.
We will use the following lemma from \cite{Carneiro}:

\begin{Lem}\label{Carneiro radial}
If $L$ is an autonomous Tonelli Lagrangian on a closed manifold $M$, then at every $h \in \rH_1(M;\R) \setminus \{0\}$, $\beta$ is differentiable in the radial direction, that is, the map 
$$
\begin{array}{rcl}
B_h: \R & \longrightarrow & \R \\
t & \longmapsto & \beta (th)
\end{array}
$$
is $C^1$ on $\R\setminus\{0\}$.
\end{Lem}

\begin{Lem}\label{Klein}
If $L$ is an autonomous Tonelli Lagrangian on the Klein bottle, then $\beta$ is $C^1$, except possibly at $0$.
\end{Lem}

\begin{Proof}
The first Betti number of the Klein bottle is one, that is, there exists $h_0 \in \rH_1(\K ;\R) \setminus \{0\}$ such that for all $h \in \rH_1(\K ;\R) $, there exists $t \in \R$ such that $h=th_0$. Then, we use Lemma \ref{Carneiro radial}.
\end{Proof}

The meaning of the next theorem is that for an autonomous Lagrangian on a closed surface, at a 1-irrational, non-singular homology class, $\beta$ is differentiable only in the directions where it is flat, and in the radial direction. Indeed, in the statement below, $\mathcal{V}(h_0)$ may be viewed as a measure of the non-differentiability of $\beta$ at $h_0$, while $\partial \alpha (c_0)$ is the largest flat containing $h_0$ in its relative interior. 
\begin{Teo}\label{betadiff, thm}
Let
\begin{itemize}
	\item $M$ be a closed surface
	\item $L$ be an autonomous  Tonelli Lagrangian on $M$
	\item $h_0$ be a 1-irrational, non-singular homology class
	\item $(\gamma_i, \dot\gamma_i)_{i \in I}$ be the periodic orbits which comprise the supports of all $(L,th)$-minimizing measures, for all $th$ in $R_{h_0}$
	\item $c_0$ be a cohomology class in the relative interior of $\partial \beta (h_0)$
	\item $\mathcal{V}(h_0)$ be the vector subspace of $\rH^1(M;\R)$ generated by the differences $c_1 - c_2$, where $c_1,c_2$ are elements of 
	$\partial \beta (h_0)$.
\end{itemize}
Then we have:
\begin{itemize}
  \item either 
  $$
  \mathcal{V}(h_0) = \partial \alpha (c_0)^{\perp}= \bigcap_{i \in J } h_i ^{\perp}
  $$
  where orthogonality is meant with respect to the duality between $\rH_1(M;\R)$ and $\rH^1(M;\R)$
  \item or $M= \T^2$   and  the closed curves $\gamma_i $, $i \in I$, foliate $M$; in this case $ \mathcal{V}(h_0) = \{0\}$. 
\end{itemize}
\end{Teo}

\begin{Proof}
We begin by proving that we always have 
 \begin{equation}\label{betadiff, 1}
  \mathcal{V}(h_0) \subseteq \partial \alpha (c_0)^{\perp},
  \end{equation}
regardless of the irrationality or singularity of $h_0$.  Take $c$ in $\partial \beta (h_0)$ and $h$ in $\partial \alpha (c_0)$. Recall that by \cite[Lemma A.1]{Massart09}, since $c_0$ lies in the relative interior of $\partial \beta (h_0)$, and $c \in \partial \beta (h_0)$, we have 
$\partial \alpha (c_0) \subseteq \partial \alpha (c)$, whence $h \in \partial \alpha (c)$. Therefore we have
\begin{eqnarray*}
\alpha(c_0)+\beta(h) & = & \langle c_0,h \rangle \\
\alpha(c)+\beta(h) & = & \langle c,h \rangle .
\end{eqnarray*}
On the other hand, since $h \in \partial \alpha (c_0)$, by \cite{Carneiro}, $\alpha(c_0)$ is the energy level of any $(L,h)$-minimizing measure. Since we also have $h \in \partial \alpha (c)$, it follows that $\alpha(c)=\alpha(c_0)$, whence 
$$
 \langle c_0 - c,h \rangle =0
 $$
 that is, $c_0-c$ lies in $h^{\perp}$. Since $h$ is arbitrary in $\partial \alpha (c_0)$, this yields $c_0-c \in \partial \alpha (c_0)^{\perp}$, whence (\ref{betadiff, 1}).
 
 Now we prove that under the hypothesis of the theorem, we have 
  \begin{equation}\label{betadiff, 2}
  \mathcal{V}(h_0) \supseteq \partial \alpha (c_0)^{\perp},
  \end{equation}
  unless $M=\T^2$ and the $\gamma_i$ foliate $\T^2$. 
  First observe that by Lemma \ref{h_i conv partial alpha}, 
  $$
  \partial \alpha (c_0)^{\perp} = \bigcap_{i \in J } h_i ^{\perp}.
  $$
  Since $c_0$ lies in the relative interior of  $\partial \beta (h_0)$,  \cite[Lemma A.4]{Massart09} says that the largest flat of $\alpha$ containing $c_0$ in its relative interior is $\partial \beta (h_0)$. \\

Recall  \cite[Theorem 1]{Massart03}:
{\it Let $L$ be an autonomous Tonelli Lagrangian on a compact surface $M$, for any cohomology class $c$. If
\begin{itemize}
  \item $F_c$ is the largest flat of $\alpha$ containing $c$ in its relative interior
  \item $V_c$ is the vector space generated by differences $c'-c''$, for all $c', c''$ in $F_c$
  \item $E_c$ is the space of cohomology classes of $1$-forms supported away from $\cA_c$,
\end{itemize}
then 
$$
E_c = V_c.
$$
}

\vspace{10 pt}

In our case, $F_{c_0} = \partial \beta (h_0)$, so $V_{c_0} =  \mathcal{V}(h_0)  $, hence 
\begin{equation}\label{betadiff, 4}
E_{c_0}=  \mathcal{V}(h_0) . 
\end{equation}
Now we will prove that either the $\gamma_i$ foliate $M$, or $E_{c_0} = \bigcap_{i \in J } h_i ^{\perp}$.
First recall that by Theorem \ref{AvsM, principal}, the projected Aubry set $\mathcal{A}_{c_0}$  is precisely the union of the $\gamma_i$.

Denote by $A_{\epsilon}$ the $\epsilon$-neighborhood of $\cA_{c_0}= \cup_{i \in I} \gamma_i$ in $M$. 
Since $\cA_{c_0}$ consists of pairwise disjoint, simple closed curves, either \begin{itemize}
  \item $\cA_{c_0}=M$, 
  \item or there exists $\epsilon_1 >0$ such that for any $0< \epsilon \leq \epsilon_1 $, $A_{\epsilon}$ is a disjoint union of annuli and M\"{o}bius strips.
  \end{itemize} 
  In the first case, $M$ has a non-singular foliation by closed curves, so $M$ is the torus or the Klein bottle. If $M$ is the Klein bottle, and the 
  $\gamma_i$ foliate $M$, then the homology classes $h_i$ are non-zero and pairwise proportional. Since the first Betti number of the Klein bottle is one, we have $\bigcap_{i \in J } h_i ^{\perp}= \{0\}$. On the other hand, by Lemma \ref{Klein},  $\beta$ is $C^1$ so $\mathcal{V}(h_0)= \{0\}$.
  Thus the conclusion of the theorem holds when $\cA_{c_0}=M$. 
  
  So assume we are in the second case. 
Observe that since the dimension of $\rH^1 (M;\R)$ is finite, there exists $\epsilon_0 >0$ such that for any $0< \epsilon \leq \epsilon_0 $, $E_c$ is the set of cohomology classes of $1$-forms supported at least $\epsilon$-away from $\cA_c$. Take $0<\epsilon \leq  \epsilon_0$. Set 
$$
\cH_\e:= i_* (\rH_1 (A_{\epsilon}; \R)),
$$
where $i :\   A_{\epsilon} \longrightarrow M$ is the canonical inclusion. Then $\cH_\e$ is generated by the homology classes $h_i$, where $i$ lies in the finite set $J$. Thus, denoting
$$
\cH_\e^{\perp}:= \{ c \in \rH^1(M;\R) :\   \forall h \in \cH_\e, \  \langle c,h \rangle =0 \}, 
$$
we have 
\begin{equation}\label{betadiff, 3}
\cH_\e^{\perp}=\bigcap_{i \in J } h_i ^{\perp}.
\end{equation}
So, in order to prove the theorem, we just need to show that $\cH_\e^{\perp}=E_{c_0}$. Take $c$ in $\cH_\e^{\perp}$ and a closed $1$-form $\eta$ with cohomology $c$. Then $\eta$ vanishes on any cycle contained in $A_{\epsilon}$, so there exists a $C^1$ function $f$ on $M$ such that $\eta_x  = d_x f$ for every $x$ in $A_{\epsilon}$. Then, $\eta -df $ is supported outside $A_{\epsilon}$, so $\left[ \eta -df \right]=c$ lies in $E_{c_0}$. This proves 
$\cH_\e^{\perp}\subseteq E_{c_0}$. The other inclusion is obvious, so $\cH_\e^{\perp}=E_{c_0}$.  By Equations (\ref{betadiff, 4}) and  (\ref{betadiff, 3}), this yields Equation (\ref{betadiff, 2}), and finishes the proof of the theorem. 

\end{Proof}

From Theorem \ref{betadiff, thm} we deduce that when $M$ is oriented, $\beta$ is never differentiable at any 1-irrational, non-singular homology class, unless $M= \T^2$ and $M$ is foliated by periodic orbits.
\begin{Cor}\label{oriented, not differentiable}
Let
\begin{itemize}
	\item $M$ be a closed, oriented  surface
	\item $L$ be an autonomous  Tonelli Lagrangian on $M$
	\item $h_0$ be a 1-irrational, non-singular homology class
	\item $(\gamma_i, \dot\gamma_i)_{i \in I}$ be the periodic orbits which comprise the supports of all $(L,th)$-minimizing measures, for all $th$ in $R_{h_0}$.
\end{itemize}
Then we have :
\begin{itemize}
  \item 
  either the dimension of $\partial \beta (h_0)$ is at least the genus of $M$
    \item or $M= \T^2$ and  the closed curves $\gamma_i $, $i \in I$, foliate $M$; in this case $\beta$ is differentiable at $h_0$. 
\end{itemize}
\end{Cor}
\begin{Proof}
We use the notations of the proof of Theorem \ref{betadiff, thm}.
 Recall that, since $M$ is oriented,  the first homology of $M$ is  endowed with a symplectic form, induced by the algebraic intersection of oriented cycles. Since $\cH_\e$ is generated by the homology classes of  disjoint closed curves, $\cH_\e$ is isotropic with respect to the symplectic intersection form, hence the dimension of $\cH_\e$ is at most half that of $\rH_1(M;\R)$. Recall that  half the dimension  of $\rH_1(M;\R)$ is precisely the genus of $M$. So the dimension of $\cH_\e^{\perp}=E_c$ is at least the genus of $M$. Now, we have seen in the proof of Theorem \ref{betadiff, thm} that either the $\gamma_i $, $i \in I$, foliate $M$, or $\cH_\e^{\perp}$ is precisely the vector subspace of $\rH^1(M;\R)$ generated by the differences of elements of $\partial \beta (h_0)$, thus proving that the dimension of $\partial \beta (h_0)$ is at least the genus of $M$.
 \end{Proof}
 
 \begin{Rem} 
Observe that the above result is not true if $h$ is singular. Indeed, take a vector field $X$ on a closed surface $\Sigma_2$ of genus $2$, which consists of periodic orbits, except for a singular graph with two fixed points and four hetero/homoclinic orbits between the two fixed points (see figure 1).
One can embed it into the Euler-Lagrange flow of a Tonelli Lagrangian, given by 
$L_X(x,v)=\frac{1}{2}\|v-X(x)\|^2_x$. It is possible to show (see for instance \cite{FFR}) that ${\rm Graph}(X)$ is invariant and that in this case $\widetilde{\cA}_0={\rm Graph}(X)$; therefore, the projected Aubry set $\cA_0$ is the whole surface $\Sigma_2$,
so  $\beta$ is differentiable at all homology classes in $\dpr \a (0)$ (see \cite[Theorem 1]{Massart03}). However,  $\cA_0=\Sigma_2$ is not foliated by periodic orbits.
\begin{figure} [h!]
\begin{center}
\includegraphics[scale=3.5]{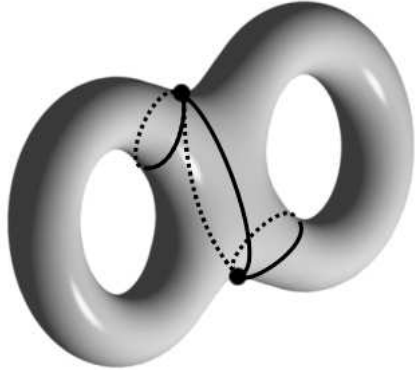}
\caption{A surface $\Sigma_2$ of genus two, consisting of periodic orbits, except for  two fixed points and four heteroclinics between them.}
\end{center}
\end{figure}
\end{Rem}

 When $M$ is not orientable, the situation is different because $\beta$ may have flats of maximal dimension (that is, of dimension equal to the first Betti number of $M$). So $\beta$ may well be differentiable at some 1-irrational, non-singular homology class. 
 For instance, by \cite[Theorem 1.3]{nonor},  there exists a Riemannian metric on $M$, whose stable norm has a (finite) polyhedron as its unit ball. Let $L$ be the Lagrangian induced by this Riemannian metric, then the $\beta$-function of $L$ is half the square of the stable norm. In particular $\beta$ is differentiable everywhere but along a finite number of straight lines. So $\beta$ is differentiable at most 1-irrational classes. Furthermore, since the Lagrangian is a Riemannian metric:
\begin{itemize}
  \item no homology class other than zero is singular
  \item the $\beta$-function is quadratic (\ie 2-homogeneous) so radial faces are trivial, \ie $\forall h, \  R_h = \{h\}$. 
\end{itemize} 
On the other hand, for every homology class the Mather set is a finite union of closed geodesics, so we never have $\cM (R_h)=M$. \\
In general all we have is the following.

  \begin{Cor}\label{maintheo2, nonor}
Let $M$ be a closed, non-orientable surface other than the projective plane or the Klein bottle,  and let $L:\rT M \longrightarrow \R$ be a Tonelli Lagrangian. Then, there exists some homology class $h$ such that $\beta$ is not differentiable at $h$.
\end{Cor}

\begin{Proof}
Let  
\begin{itemize}
  \item  $h_0\in \rH_1(M;\R)$ be 1-irrational and   non-singular
  \item $c_0$ be a cohomology class in the relative interior of $\partial \beta (h_0)$
  \item $(\gamma_i, \dot\gamma_i)_{i \in I}$ be the periodic orbits which comprise the supports of all $(L,th)$-minimizing measures, for all $th_0$ in $R_{h_0}$
     \item $\mu_i$ be the probability measure equidistributed on $(\gamma_i, \dot\gamma_i)$
  \item $h_i$ be the homology class of $\mu_i$
  \item $J$ be the finite subset of $I$, given by Lemma \ref{h_i conv partial alpha},  such that $\partial \alpha (c_0)$ is the convex hull of the $h_i$, $i \in J$.
\end{itemize}
We are going to show that $\beta$ is not differentiable at $h_j$, for some $j \in J$. Because of Theorem  \ref{betadiff, thm} we need to choose $h_j$ carefully, so that $R_{h_j}$ does not meet the relative interior of $\partial \alpha (c_0)$. This can be done as follows. 
By the Hahn-Banach theorem, there exists an affine hyperplane $\cH$ in $\rH_1(M;\R)$, which separates  $\partial \alpha (c_0)$ from zero. Define 
$$
\begin{array}{rcl}
\Psi :\  \partial \alpha (c_0) & \longrightarrow & \cH \\
h & \longmapsto & \R h \cap \cH.
\end{array}
$$
Then $\Psi (\partial \alpha (c_0))$ is convex because $\Psi$ sends straight lines to straight lines. Take an extremal point $x$ in 
$\Psi (\partial \alpha (c_0))$. Then for all $h \in \Psi^{-1}(x)$, $h$ is an extremal point of $\partial \alpha (c_0)$, so there exists some $j$ in $J$ such that $\Psi (h_j)= x$.

First note that $h_j$ is 1-irrational because $h_j = \left[ \mu_j \right] = T_j ^{-1} \left[ \gamma_j \right]$, where $T_j$ is the least period of the periodic orbit $(\gamma_j, \dot\gamma_j)$. 
Furthermore $h_j$ is non-singular, as we now prove. 
Take $c$ in the relative interior of $\partial \beta (h_j)$. Since $h_j \in  \partial \alpha (c_0)$, we have $c_0 \in \partial \beta (h_j)$. Now by \cite{Massart03}, since $c_0 \in \partial \beta (h_j)$ and $c$ lies in the relative interior of $\partial \beta (h_j)$, the Aubry set of $c$ is contained in the Aubry set of $c_0$, which is the union of the $(\gamma_i, \dot\gamma_i)$ over al $i \in I$ by Theorem  \ref{AvsM, principal}. Therefore the Aubry set of $c$ does not contain any fixed point. Since  $c$ lies in the relative interior of $\partial \beta (h_j)$, this entails that $h_j$ is non-singular. 

We want to apply Theorem \ref{betadiff, thm} to $h_j$. Take $t$ such that $th_j$ lies in $R_{h_j}$. By \cite[Lemma 2.3]{Massart09}, we have 
$\partial \beta (h_j) \subseteq \partial \beta (th_j) $, so $c_0 \in \partial \beta (th_j) $, that is, $th_j \in  \partial \alpha (c_0)$. Now recall that $\partial \alpha (c_0)$ is the convex hull of the $h_i, i \in I$, so any $(L,th_j)$-minimizing measure is supported on some of the  $(\gamma_i, \dot\gamma_i)$. Let $K$ be the subset of $I$ such that the union of the supports of all $(L,th_j)$-minimizing measures, for all $t$ such that $th_j$ lies in $R_{h_j}$, is the union over $ k \in K$ of $(\gamma_k, \dot\gamma_k)$. Now we prove that for all $k \in K$, $h_k \in \R h_j$. 

Assume it is not so, \ie there exists some $k$ in $K$ such that 
$h_k = \left[\mu_k\right] \not\in \R h_j$. Since $k \in K$, there exists $t$ such that $th_j \in R_{h_j}$, and $\mu_k$ is an ergodic component of some $(L,th_j)$-minimizing measure. Thus we may write 
$$
th_j = \lambda h_k + (1-\lambda) h,
$$
where $\lambda \in \left( 0,1 \right)$ and $h \in \partial \alpha (c_0)$. Then $\Psi (th_j)$, which is $x$, is a convex combination of $\Psi (h_k)$, which is not $x$, and $\Psi (h)$. This contradicts the extremality of $x$. We have proved that for all $k \in K$, $h_k \in \R h_j$. Now we apply Theorem \ref{betadiff, thm} to $h_j$:  since $M$ is not the Klein bottle, we have  
$$
  \mathcal{V}(h_j) =  \bigcap_{k \in K} h_k ^{\perp}.
  $$
Now since for all $k \in K$, $h_k \in \R h_j$, we have  
$$
 \bigcap_{k \in K} h_k ^{\perp}= h_j ^{\perp}.
 $$
  Since $M$ is neither the Klein bottle nor the projective plane,  the dimension of $h_j ^{\perp}$ is greater than zero, that is, $\beta$ is not differentiable at $h_j$.
 \end{Proof}

For the sake of completeness, let us recall that by \cite[Corollary 3]{Massart03} we also have:

\begin{Prop}\label{prop2}
For any autonomous Lagrangian $L$ on a closed surface $M$, 
$\b$ is  differentiable at every completely irrational homology class.
\end{Prop}

 %%%%%%%%%%%%%%%%
 \subsection{The case when $\beta$ is $C^1$ and $M=\T^2$}
  When $M$ is the two-torus and  $\beta$ is $C^1$ everywhere, we can rule out radial flats of $\beta$.

 \begin{Prop} \label{C1_implies_convex}
Let $L$ be an autonomous Tonelli Lagrangian on the two-torus such that $\beta$ is differentiable at every point of $\rH_1(\T^2;\R)$. Then for all $h \in \rH_1(\T^2;\R) \setminus \{ 0 \}$, we have $R_h = \{ h \}$. As a consequence, $\beta$ is strictly convex.
\end{Prop}
\begin{Proof} 

{\bf First case:} $h$ is 1-irrational.

Replacing, if we have to,  $h$ with an extremal point of $R_h$, which is also 1-irrational, we may assume that $R_h = \left[ \lambda h, h \right]$ for some $0 \leq \lambda \leq 1$. We want to prove that $\lambda =1$. Pick \begin{itemize}
  \item $(x,v) \in \rT\T^2$ such that the probability measure carried by the orbit $\Phi^L_t (x,v) $ has homology $\lambda h$
    \item a sequence of real numbers $t_n > 1$ such that $t_n$ converges to $1$
  \item for each $n \in \N$,  $c_n \in \rH^1 (\T^2;\R)$ such that $\partial \beta (t_nh) =c_n$ (recall that $\beta$ is differentiable at $t_nh$).
  \end{itemize}
  First note that $t_n h $ is not singular, for any $n$. Indeed if it were, then zero would lie in the radial flat of $t_n h$, that is, $\beta$ would be affine along the segment $\left[ 0, t_n h\right]$. But $t_n >1$, so  the segment $\left[ 0, t_n h\right]$ properly contains $\left[ \lambda h, h \right]$, which  contradicts  $R_h = \left[ \lambda h, h \right]$.
Then by Corollary \ref{oriented, not differentiable}, since $\beta$ is differentiable at $t_nh$, which is 1-irrational, we have $\cM_ {c_n}=M$, thus there exists, for every $n$, a $v_n \in T_x M$ such that $(x,v_n) \in  \widetilde{\cM}_{c_n} $, and moreover the orbit $\Phi^L_t (x,v_n) $ is periodic. By semicontinuity of the Aubry set, and by the Graph Theorem, $v_n$ converges to $v$ as $n$ goes to infinity. Let $T$ and $T_n$ be the minimal periods of 
$\Phi^L_t (x,v) $ and $\Phi^L_t (x,v_n) $, respectively. If $(x,v)$ is a fixed point we just set $T := + \infty$. We now prove that $\liminf T_n \geq T$.

Indeed, if  some subsequence $T_{n_k}:= T_k$ converged to $T' < T$, we would have 
$\Phi^L_{T' }(x,v) = (x,v)$, contradicting the minimality of $T$. If $(x,v)$ is a fixed point, we have $v=0$, so $v_n$ converges to zero, so $T_n$ tends to infinity.

Let $h_0 \in \rH_1(\T^2;\Z)$ be a generator of $\R h \cap \rH_1(\T^2;\Z)$ such that the probability measures carried by the orbits $\Phi^L_t (x,v) $ and $\Phi^L_t (x,v_n) $ have homology $T^{-1} h_0 = \lambda h$ and $T_n ^{-1} h_0$, respectively. 

Now  if 
$\lambda$ were $< 1$, since $\liminf T_n \geq T$, for $n$ large enough there would exist a   $c_n$-minimizing measure with homology in $\left[ \lambda h, h \right)$. This means that the radial flats $R_h$ and $R_{t_n h}$ overlap, in other words, 
$t_n h \in R_h$. This contradicts the fact that $R_h = \left[ \lambda h, h \right]$, hence 
$\lambda =1$.

{\bf Second case:} $h$ is 2-irrational, that is, completely irrational. Then any $(L,h)$-minimizing measure is supported on a lamination of the torus without closed leaves. Any such lamination is uniquely ergodic, in particular $h$ is not contained in any non-trivial flat of $\beta$, radial or not,  regardless of the Lagrangian.

The statement about the differentiability of $\beta$ implying its strict convexity is now just a consequence of the fact, observed by Carneiro in \cite{Carneiro}, that for an autonomous Lagrangian on a closed, orientable surface $M$, the flats of $\beta$ are contained in isotropic subspaces of $\rH_1 (M;\R)$ with respect to the intersection symplectic form on $\rH_1 (M;\R)$. In particular, when $M = \T^2$, all flats are radial.
 \end{Proof}
 
 \begin{Prop}\label{lemma0}
Let  $L: \rT \T^2 \longrightarrow \R$ be a Tonelli Lagrangian. Assume that   $\beta$ is $C^1$ everywhere.
Then zero is the only (possibly) singular class, and for every 1-irrational, nonzero homology class $h$, $\cM^h = \T^2$, and $\cM^h$ is foliated by periodic orbits with the same homology class and minimal period. 
\end{Prop}

\begin{Proof}
First note that if a non-zero class $h$ is singular, then there exists $c$ in $\partial \beta (h)$, such that there is a fixed point in the Mather set of $c$. Thus $R_h$ contains the homology of the Dirac measure on the fixed point, which is zero. This contradicts Proposition \ref{C1_implies_convex}, which says that $R_h = \left\{ h \right\}$.

Now take a non-zero, 1-irrational homology class $h$, so $h$ is non-singular, and by Corollary \ref{oriented, not differentiable}, $\cM (R_h)=\T^2$, and $\cM (R_h)$ is foliated by periodic orbits. Since  $R_h = \left\{ h \right\}$ by Proposition \ref{C1_implies_convex}, we have $\cM^h = \T^2$. It is well known that all leaves of a foliation of the two-torus by homologically non trivial closed curves are homologous. So there exists $h_0 \in \rH_1 (\T^2;\Z)$ such that the projection to $\T^2$ of any orbit contained in $\cM^h$ is homologous to $h_0$. Besides, for each $x \in \T^2$, let $T(x)$ be the minimal period of the periodic orbit $\gamma_x$ in $\cM^h$ through $x$. The homology class of the probability measure carried by $\gamma_x$ is 
$T^{-1}h_0$. Now the fact that there are no non-trivial radial flats of $\beta$ implies that $T(x)$ is independent of $x$.
\end{Proof}

%%%%%%%%%%%%%%%%%%%%%%%%%%%%%%%%%%

\section{Integrability}\label{section4}

Let us  recall the definition of $C^0$-integrability (see also \cite{Arnaud}).

\begin{Def}\label{C0integrability}
A Tonelli Hamiltonian $H:\rT^*M \longrightarrow \R$ is said to be $C^0$-integrable, if there exists a foliation of $\rT^* M$ made by invariant Lipschitz Lagrangian graphs, one for each cohomology class.
\end{Def}

\begin{Lem}\label{lemma2}
Let $M$ be a compact manifold of any dimension, $L: \rT M \longrightarrow \R$ a Tonelli Lagrangian and $H: \rT^* M \longrightarrow \R$ the associated Hamiltonian. If $H$ is $C^0$-integrable, then $\beta$ is $C^1$.
\end{Lem}

\begin{Proof}
Suppose that $H$ is $C^0$-integrable. This means that the cotangent space $\rT^*M$ is foliated by disjoint invariant Lipschitz Lagrangian graphs. Let us denote  by $\L_c$ the invariant Lagrangian graph of such a foliation, corresponding to the cohomology class $c$. It is easy to check that each $\L_c$ is the graph of a solution of Hamilton-Jacobi equation $H(x,\eta_c+du)=\a(c)$, where $\eta_c$ is a closed $1$-form on $M$ with cohomology class $c$,
 and therefore
from weak KAM theory \cite{Fathibook} it follows that for each $c\in \rH^1(M;\R)$ the Aubry set $\cA^*_c \subseteq \L_c$.
If for some $h\in\rH_1(M;\R)$ there exist $c\neq c'$ such that $c,c'\in\dpr \b(h)$, then  $\widetilde{\cM}^h \subseteq \widetilde{\cM}_c\cap \widetilde{\cM}_{c'}$ (Lemma \ref{lemma1}). But this implies that $\widetilde{\cM}_c\cap \widetilde{\cM}_{c'} \neq 0$ and therefore 
$
\L_c \cap \L_{c'} \supseteq \cA_c^* \cap \cA^*_{c'} \supseteq \cM_c^* \cap \cM^*_{c'} \neq 0,
$
which contradicts the disjointness of the Lagrangian graphs forming the foliation. Hence, $\b$ must be differentiable everywhere.\\
\end{Proof}

We can now prove the main results stated in Section \ref{section1}.\\

\begin{Prop} \label{torus_only_integrable}
The torus is the only closed surface which admits a $C^0$-integrable Hamiltonian.
\end{Prop}
\begin{Proof}

First,  no Hamiltonian on the sphere can be $C^0$-integrable. Indeed, any Lagrangian graph is exact since the sphere is simply connected, and any two exact Lagrangian graphs intersect, because any $C^1$ function on the sphere has critical points. Likewise, no Hamiltonian on the projective plane can be 
$C^0$-integrable, otherwise its lift to the sphere would be $C^0$-integrable.

For the Klein bottle $\K$, we need to work a little bit more. 
 For each $x \in \K$, let us define
\beqano
F_{x} : \rH^1(\K;\R)\simeq \R &\longrightarrow & \rT^*_x \K \simeq \R^2\\
c &\longmapsto& \L_c \cap \rT^*_x \K,
\eeqano
where $\L_c$, for  $c\in\rH^1(\K;\R)$, are the Lagrangian graphs foliating $\rT^*\K$.
It is possible to check that $F_{x}$ is continuous (see \cite[Lemme 4.22]{Arnaud} or the proof of the continuity of (\ref{function}) in Theorem \ref{maintheo1}) and  injective (as it follows from the disjointness of the $\L_c$'s). Moreover, if the Hamiltonian is $C^0$-integrable, the map $F_x$ is surjective. Now there is no such thing as a continuous bijection from $\R$ to $\R^2$, so there is no $C^0$-integrable Hamiltonian on the Klein bottle. 

The same argument can be used for any surface with first Betti number $>2$. Finally,
no Hamiltonian on the connected sum of three projective planes can be $C^0$-integrable, otherwise it would  lift to  a $C^0$-integrable Hamiltonian on a surface of genus two.

\end{Proof}

\begin{Teo}\label{maintheo1}
Let  $L:\rT \T^2 \longrightarrow \R$ be a Tonelli Lagrangian on the two-torus. 
Then, $\beta$ is $C^1$ if and only if the system is $C^0$-integrable.
\end{Teo}

\begin{Proof}
[$\Longleftarrow$] It follows from Lemma \ref{lemma2}. 
 [$\Longrightarrow$] For each homology class $h$, let us denote by $c_h:=\dpr \beta(h)$. 
% Observe that since $\beta$ is differentiable everywhere, then if $h\neq 0$ then it is non-singular. 
 If $h$ is non-singular and $1$-irrational, then  Proposition \ref{lemma0} says that $\L_{c_h}:= \cA^*_{c_h}$ is an invariant Lipschitz Lagrangian graph of cohomology class $c_h$, which is foliated by periodic orbits of homology $h$ and same minimal period $T_h$. Since the dynamics on these Lipschitz Lagrangian graph is totally periodic, \ie $\Phi^H_{T_h}\big| \L_{c_h} = {\rm Id}\big| \L_{c_h}$, then it follows from the result in \cite{Arnaud} (see for instance the proof of Th\'eor\`eme 4) that this graph is in fact $C^1$.
 
Observe that such cohomology classes $c_h$ are dense in $\rH^1(\T^2;\R)$. Indeed, 1-irrational non-singular homology classes are dense in $\rH_1(\T^2;\R)$ 
 because 1-irrational  homology classes are dense in $\rH_1(\T^2;\R)$ and, by Proposition \ref{lemma0}, zero is the only singular class.
 On the other hand, since $\beta$ is differentiable, the Legendre transform is a homeomorphism from $\rH_1(\T^2;\R)$ to $\rH^1(\T^2;\R)$.

   Using the semicontinuity of the Aubry set in dimension $2$ (see \cite{Bernard08}), we can deduce that for each $c\in \rH^1(\T^2;\R)$ the Aubry set $\cA^*_c$ projects over the whole manifold and therefore it is still an invariant Lipschitz Lagrangian graph, which we will denote $\L_c$. Observe that all these $\L_c$'s are disjoint. This is a straightforward consequence of the differentiability of $\beta$. In fact if for some $c\neq c'$ we had that $\L_c \cap \L_{c'} \neq 0$, then  ${\cM}^*_c\cap {\cM}^*_{c'} \neq 0$; but this would contradict the differentiability of $\beta$ at some homology class. It only remains to prove that the  union is the whole $\rT^*\T^2$. Observe that for each $c\in \rH^1(\T^2;\R)$, $\L_c$ is the graph of the unique solution of the Hamilton-Jacobi equation $H(x,\eta_c + du_c) = \a(c)$, where $u_c\in C^{1,1}(\T^2)$ and $\eta_c$ is a closed $1$-form on $\T^2$ with cohomology class $c$.  For each ${x_0}\in \T^2$, let us define
\beqa{function}
F_{x_0} : \rH^1(\T^2;\R)\simeq \R^2 &\longrightarrow & \rT^*_{x_0}\T^2 \simeq \R^2\\
c &\longmapsto& \eta_c(x_0) + d_{x_0}u_c.\nonumber
\eeqa
This map is injective (as it follows from the disjointness of the $\L_c$'s). Moreover,  $F_{x_0}$ is also continuous. In fact,  let $c_n\rightarrow c$ and consider the sequence of Lipschitz functions $\lambda_{c_n}:= \eta_{c_n} + du_{c_n}$. This sequence  is equibounded and equilipschitz (it follows from the results in \cite{Fathibook} or from Mather's graph theorem for the Aubry sets \cite{Mather93}). Therefore, applying Ascoli-Arzel\`a's theorem, we can conclude that - up to selecting a subsequence - $\lambda_{c_n}$ converge uniformly to some $\tilde{\lambda}= \eta_c+ du$, for some $u\in C^{1,1}(\T^2)$. Observe that since $H(x,\lambda_{c_n}(x))=\a(c_n)$ for all $x\in \T^2$ and all $n$, and $\a$ is continuous, then $H(x,\tilde{\lambda}(x))=\a(c)$ for all $x$. Therefore, $u$ is a solution of Hamilton-Jacobi equation $H(x, \eta_{c}+du)=\a(c)$. But, as we noticed above, for any given $c\in \rH^1(\T^2;\R)$ there exists a unique solution of this equation (essentially because $\cA_c=\T^2$). Therefore, $\tilde{\l}=\lambda_c$ and hence $F_{x_0}(c_n) \longrightarrow F_{x_0}(c)$.

Using that the above map is continuous, we can conclude that $F_{x_0}(\R^2)$ is open (see for instance \cite{Brouwer}). In a similar way, one can show that this image is also closed. In fact, let $y_k = F_{x_0}(c_k)$ be a sequence in $\rT^*_x \T^2$ converging to some $y_0$. Since each family of classical solutions of Hamilton-Jacobi on which the $\a$-function (\ie the energy) is bounded gives rise to a family of equi-$C^{1,1}$ functions (see \cite{Fathibook}),  then there exist
$\tilde{u} \in C^{1,1}(\T^2)$ and $\tilde{c} \in \rH^1(\T^2;\R)$, such that 
$c_k+du_{c_k} \rightarrow \tilde{c}+d\tilde{u}$ uniformly. From the continuity of the $\a$-function, it follows that $H(x,\tilde{c}+d_x\tilde{u})=\a(\tilde{c})$. As before, since $\L_{\tilde{c}}$ is the unique invariant Lagrangian graph with cohomology $\tilde{c}$, then $\L_{\tilde{c}}=\{(x,\tilde{c}+d_x\tilde{u}:\; x\in \T^2\}$ and therefore $y_0 = \tilde{c} + d_{x_0}\tilde{u} = F_{x_0}(\tilde{c}) \in F_{x_0}(\R^2)$. This shows that $F_{x_0}(\R^2)$ is also closed and therefore it is all of $\R^2$. Since this holds for all $x_0\in \T^2$, we can conclude that $\bigcup_c \L_c = \rT^*\T^2$, that is, the system is $C^0$-integrable. 
\end{Proof}

We can deduce more information about the dynamics on a $C^0$-integrable system (see also \cite{Arnaud}).

\begin{Cor}\label{corollarymaintheo}
Let $H: \rT^* M \longrightarrow \R$ be a $C^0$-integrable Hamiltonian on a two-dimensional closed manifold $M$. Then, $M$ is diffeomorphic to $\T^2$. Moreover :
\begin{itemize}
\item[{\rm (}i{\rm)}]  for each $1$-irrational homology class $h$, there exists an invariant Lagrangian graph foliated by periodic orbits with homology $h$ and the same minimal period;
\item[{\rm (}ii{\rm)}] for each completely irrational homology class, there exists an invariant Lagrangian graph on which the motion is conjugated to an irrational rotation on the torus or to a Denjoy type homeomorphism; 
\item[{\rm (}iii{\rm)}] there exists a dense $G_{\delta}$ set of (co)homology classes, for which the motion on the corresponding invariant torus is conjugated to a rotation;
\item[{\rm (}iv{\rm)}] as for the $0$-homology class, there exists a $C^1$ invariant torus  $\L_{c(0)}= \{(x,\frac{\dpr L}{\dpr v}(x,0):\; x\in \T^2\}$ consisting of fixed points.
\end{itemize}
\end{Cor}

\begin{Proof}
If $H$ is $C^0$ integrable, then $\b$ is $C^1$ (Lemma \ref{lemma2}). By Proposition \ref{torus_only_integrable}, $M$ is diffeomorphic to the 2-torus and it follows from Theorem \ref{maintheo1} that for each $c\in \rH^1(\T^2;\R)$ the Aubry set $\cA^*_c$ is a Lipschitz Lagrangian graph and all these sets provide a partition of the cotangent space $\rT^*M$.

%Observe that $M$ cannot be non-orientable. In fact, if $M$ were the connected sum of three projective planes, then there could not be a foliation of $M$ by closed curves, because it would lift to a foliation of the genus $2$ orientation cover by closed curves.\\

\noindent
({\it i})  It is a consequence of Proposition \ref{lemma0}.
\noindent ({\it ii}) Let $h$ be a completely irrational homology class and let $c_h=\dpr \beta(h)$. We will follow the discussions in \cite[p.~1084]{Carneiro} (see also \cite[p.~66 et seqq.]{Mather09} for more technical details). The Hamiltonian flow on $\L_{c_h}$ induces a flow on $\T^2$ without fixed points (otherwise the homology class would be singular  and this would contradict Proposition \ref{lemma0}) or closed trajectories (otherwise the homology class would be $1$-irrational). Thus, one can construct a non-contractible closed curve $\Gamma$, which is  transverse to the flow. Using an argument {\it \`a la} Poincar\'e-Bendixon, it follows that all trajectories starting on $\Gamma$ must return to it. Hence,  one  can define a continuous map from a compact subset of $\Gamma$ to itself, which is order preserving and with irrational rotation number. Therefore, it is either conjugate to an irrational rotation or a Denjoy type homeomorphism.\\
\noindent({\it iii}) We know from ({\it i}) that for all $1$-irrational homology classes, the motion on the corresponding graph is conjugated to a rotation. Observe now that because of ({\it ii}), it is sufficient to show that the set of homology classes $\cH_r$ for which the Mather set projects over the whole torus is a dense $G_{\delta}$ set. In fact, if this is the case, then the flow cannot be conjugated to a Denjoy type homeomorphism and must be therefore conjugated to a rotation. Let us start by observing that this set $\cH_r$ is clearly dense since it contains all $1$-irrational homology classes. The fact that it is also $G_{\delta}$ follows from \cite[Corollaire 4.5]{FathiHerman}, stating  that 
 the set of strictly ergodic flows on a compact set is $G_{\delta}$ in the $C^0$ topology.
 In fact, let us consider $\cF(M)$ the set  of flows on $M$, endowed with the $C^0$ topology. One can consider the continuous and injective map given by:
 \beqano
 \psi: H_1(M;\R) &\longrightarrow& \cF(M)\\
 h &\longmapsto& \pi_M\circ \Phi^H\big| \Lambda_{c_h}\,,
 \eeqano
 \ie we associate to each $h\in H_1(M;\R)$, the projection of the Hamiltonian flow restricted to the invariant Lagrangian graph $\L_{c_h}$. Observe that injectivity easily follows from the fact that they have different rotation vectors, while continuity is a consequence of the  continuity of the Hamiltonian flow.
Now, let us define $\cH_r=\psi^{-1}\big(\cS(M) \cap \psi(H_1(M;\R))\big)$, where $\cS(M)$ denotes the set of strictly ergodic flows on $M$. Using that $\cS(M)$ is a $G_{\delta}$ in $\cF(M)$ \cite[Corollaire 4.5]{FathiHerman} and that $\psi$ is a homeomorphism on its image, we can then conclude that $\cH_r$ is also a $G_{\delta}$ set in $H_1(M;\R)$.
 
  %
% One could also observe that since our Lagrangian is $C^0$-integrable, by \cite[Th\'eor\`eme 3]{Arnaud}, there is a dense $G_{\delta}$ of cohomology classes with a $C^2$ weak KAM solution; and by a recent result of Fathi \cite{Fathi_Kyoto}, when the flow is of Denjoy type, there cannot be a $C^2$ weak KAM solution.  \\
\noindent({\it iv}) Since the union of the Aubry sets foliates all the phase space, then any point $(x,0)$ will
be contained in some Aubry set and therefore, using  \cite[Proposition 3.2]{FFR}, it follows that 
it is a fixed point of the Euler-Lagrange flow. Hence, the Dirac-measure $\d_{(x,0)}$ is invariant and action minimizing (since its support is contained in some Aubry set); clearly, such a measure  has rotation vector equal to $0$. Therefore, $\widetilde{\cM}^0 \supseteq \T^2\times \{0\}$ and from the graph property it follows that they coincide. Then, 
$\Lambda_{c(0)}:=\cA^*_{c(0)} = \left\{(x,\frac{\dpr L}{\dpr v}(x,0):\; x\in \T^2 \right\}$.\\
\end{Proof}

\begin{Rem}
Although we believe that the motions on all invariant Lagrangian tori must be conjugated to rotations, we have not been able to show more than $iii)$.
It remains an open question whether it is possible or not that a Denjoy type homeomorphism is embedded into a $C^0$-integrable Hamiltonian system.
\end{Rem}

\subsection{The case of mechanical Lagrangians}\label{mechanicalcase}
Recall that a mechanical  Lagrangian is $L(x,v) = 1/2 \,  g_x(v,v) +f(x)$, where $g$ is a Riemannian metric and $f$ is a $C^2$ function on $\T^2$. In this case we can bridge the gap between $C^0$ and $C^1$ integrability, using Burago and Ivanov's theorem on metrics without conjugate points \cite{BI}.
 
\begin{Prop}\label{prop_mechanical}
Let $L$ be a $C^0$-integrable mechanical Lagrangian on an $n$-dimensional torus. Then the potential $f$ is identically constant and the metric $g$ is flat. In particular, $L$ is $C^1$-integrable.
\end{Prop}

\begin{Prop}\label{prop_mechanical_2dim}
Let $L$ be a  mechanical Lagrangian on an $2$-dimensional torus, whose $\beta$-function is $C^1$. Then the potential $f$ is identically constant and the metric $g$ is flat. In particular, $L$ is $C^1$-integrable.
\end{Prop}

First we need the following lemma (see also \cite[Lemma 1]{Sorrentino08}), which will be needed for the proof of Proposition \ref{prop_mechanical}. Observe that one could also prove Proposition \ref{prop_mechanical_2dim} by using Corollary \ref{corollarymaintheo} ({\it iv}), which is however valid only for $n=2$.

\begin{Lem}
Let\begin{itemize}
  \item $M$ be a closed manifold of any dimension
  \item $L$ be an autonomous Tonelli Lagrangian on $M$, such that $L(x,v)=L(x,-v)$ for all $(x,v)$ in $TM$. 
\end{itemize}
Then the Aubry set $\widetilde{\cA}_0$ consists of fixed points of the Euler-Lagrange flow.
\end{Lem}
\begin{Proof}
Take $(x,v) \in \widetilde{\cA}_0$. Then by \cite{Fathibook} there exists a sequence of $C^1$ curves $\gamma_n \; : \; \left[0, T_n \right] \longrightarrow M$, such that 
\begin{itemize}
  \item $\gamma_n (0) = \gamma_n (T_n) = x$ for all $n\in \N$
  \item $T_n \longrightarrow +\infty$ as $n\longrightarrow +\infty$
   \item $\dot{\gamma}_n (0) \longrightarrow v$ and $\dot{\gamma}_n (T_n) \longrightarrow v$ as $n\longrightarrow +\infty$
   \item $\int_{0}^{T_n}\left( L(\gamma_n , \dot{\gamma}_n) + \alpha(0)\right)dt \longrightarrow 0$    as $n\longrightarrow +\infty$.   
   \end{itemize}
   Consider the sequence of curves $\delta_n \; : \; \left[0, T_n \right] \longrightarrow M$, such that $\delta_n (t) := \gamma_n (T_n -t)$ for all $n \in \N$, $t \in \left[0, T_n \right]$.
   Then since $L$ is symmetrical, $\int_{0}^{T_n}\left( L(\delta_n , \dot{\delta}_n) + \a(0) \right) dt \longrightarrow 0$    as $n\longrightarrow +\infty$, which proves that  $(x,-v) \in \widetilde{\mathcal{A}}_0$.  Therefore, by Mather's Graph Theorem, $v=0$, so $(x,v)$ lies in the zero section of $\rT M$. Now, by \cite[Proposition 3.2]{FFR}, for any Lagrangian, the intersection of $\widetilde{\cA}_0$ with    the zero section of $\rT M$ consists of fixed points of the Euler-Lagrange flow. This proves the lemma. 
   \end{Proof}

\begin{Proof}[{\bf Propositions \ref{prop_mechanical} and \ref{prop_mechanical_2dim}}]
 Now let us assume that the Lagrangian is mechanical. Then the only fixed points of the Euler-Lagrange flow are the critical points of the potential $f$, and the only minimizing fixed points are the minima of $f$. So if $f$ is not constant, the Lagrangian cannot be $C^0$-integrable. Furthermore, since the Lagrangian is $C^0$-integrable, every orbit is minimizing, in particular, there are no conjugate points. So by Burago and Ivanov's proof of the Hopf Conjecture \cite{BI}, the metric $g$ is flat. This proves Proposition \ref{prop_mechanical}.  Proposition \ref{prop_mechanical_2dim} is now just a consequence of Theorem    \ref{maintheo1}.
 \end{Proof}

%Let $k:= \dim H_1(M;\R) >2$. Suppose by contradiction that $\b$ is differentiable on an open subset $U\subseteq \rH_1(M;\R) \simeq \R^k$. Proceeding as before, one can show that for each $h \in \overline{U}$, $\L_{c_h}:= \cA^*_{c_h}$ is an invariant Lagrangian graph of cohomology class $c(h)$. The set of such cohomology classes, which we will denote by $V$, is a closed set with non-empty interior. Since $\b$ is differentiable on $U$, then these Lagrangian graphs are disjoint (same proof as in Proposition \ref{proptorus}). Let us now fix $x_0\in  M$ and consider the function $F_{x_0} : V  \longrightarrow  \rT*_xM \simeq \R^2$ defined as in Proposition \ref{proptorus}. This function is continous and injective, but this is impossible since $k > 2$.

%%%%%%%%%%%%%%%%%%%REFERENCES %%%%%%%%%%%%%%%%%%%
\def\cprime{$'$}

\vspace{1.truecm}

\noindent {\sc Daniel Massart}, 
{\it D\'epartement de Math\'ematiques,
Universit\'e Montpellier II, France.
\\
email: {\sc massart@math.univ-montp2.fr}}\\

\noindent {\sc Alfonso Sorrentino}, 
{\it Department of Pure Mathematics and Mathematical statistics, University of Cambridge, UK
\\
email: {\sc A.Sorrentino@dpmms.cam.ac.uk}}\\


\begin{thebibliography}{10}
\expandafter\ifx\csname natexlab\endcsname\relax\def\natexlab#1{#1}\fi
\expandafter\ifx\csname bibnamefont\endcsname\relax
  \def\bibnamefont#1{#1}\fi
\expandafter\ifx\csname bibfnamefont\endcsname\relax
  \def\bibfnamefont#1{#1}\fi
\expandafter\ifx\csname citenamefont\endcsname\relax
  \def\citenamefont#1{#1}\fi
\expandafter\ifx\csname url\endcsname\relax
  \def\url#1{\texttt{#1}}\fi
\expandafter\ifx\csname urlprefix\endcsname\relax\def\urlprefix{URL }\fi
\providecommand{\bibinfo}[2]{#2}
\providecommand{\eprint}[2][]{\url{#2}}


\bibitem{Arnaud}
Marie-Claude Arnaud.
\newblock Fibr\'es de Green et r\'egularit\'e des graphes $C^0$-Lagrangiens invariants par un flot de Tonelli.
\newblock {\em Ann. Henri Poincar\'e}, 9 (5): 881--926, 2008.

%\bibitem{TieredArnaud}
%Marie-Claude Arnaud.
%\newblock The tiered Aubry set for autonomous Lagrangian functions. 
%\newblock  {\rm Ann. Inst. Fourier}, 58 (5): 1733--1759, 2008.


%\bibitem{Arnoldbook}
%Vladimir~I. Arnol{\cprime}d.
%\newblock Mathematical methods of classical mechanics
%\newblock {\em Graduate Texts in Mathematics}, 60. Springer-Verlag, New York,  1989.


\bibitem{Bangert}
   Victor Bangert.
   \newblock Geodesic rays, Busemann functions and monotone twist maps.
   \newblock {\em Calc. Var. Partial Differential Equations}, 2 (1): 49--63, 1994.

\bibitem{nonor}
Florent Balacheff and Daniel Massart
\newblock Stable norms of non-orientable surfaces.  
\newblock {\em Ann. Inst. Fourier} (Grenoble)  58  (2008),  no. 4, 1337--1369.

%\bibitem{Bernard}
%Patrick Bernard.
%\newblock Symplectic aspects of Mather theory.
%\newblock {\em Duke Math. J.}, 136 (3): 401--420, 2007.

\bibitem{bernardc11}
Patrick Bernard.
\newblock Existence of $C^{1,1}$ critical subsolutions of the Hamilton-Jacobi equation on compact manifolds.
\newblock {\em Ann. Sci. \'Ecole Norm. Sup.}, 40 (3): 445--452, 2007.

\bibitem{Bernard08}
Patrick Bernard.
\newblock On the Conley decomposition of Mather sets.
\newblock {\em Rev. Mat. Iberoam.}, 26 (1): 115--132, 2010.


\bibitem{Brouwer}
Luitzen~E.~J.  Brouwer.
\newblock Zur invarianz des $n$-dimensionalen Gebiets.
\newblock {\em Math. Ann.}, 72 (1): 55--56, 1912.


\bibitem{BI}
Dmitri Burago and  Sergey Ivanov
\newblock Riemannian tori without conjugate points are flat.  
\newblock {\em Geom. Funct. Anal},  4  (1994),  no. 3, 259--269.

\bibitem{Carneiro}
Mario J. Dias Carneiro.
\newblock On minimizing measures of the action of autonomous {L}agrangians.
\newblock {\em Nonlinearity}, 8 (6): 1077--1085, 1995. 


\bibitem{Carneiro-Lopes}
Mario J. Dias Carneiro and Arthur  Lopes.
\newblock On the minimal action function of autonomous Lagrangians associated to magnetic fields.  
\newblock {\em Ann. Inst. H. Poincar\'e Anal. Non Lin�aire}  16  (1999),  no. 6, 667--690.

\bibitem{FathiHerman}
Albert Fathi and Michael~R. Herman.
\newblock Existence de diff\'eomorphismes minimaux.
\newblock  {\em Ast\'erisque} 49, 37--59, 1977. 

\bibitem{Fathibook}
Albert Fathi.
\newblock The Weak {K}{A}{M} theorem in {L}agrangian dynamics.
\newblock 10th preliminary version, 2008

%\bibitem{Fathi_Kyoto}
%Albert Fathi.
%\newblock {D}enjoy-{S}chwartz and {H}amilton-{J}acobi.
%\newblock {\em RIMS meeting: Viscosity solutions of differential equations and related topics}, Kyoto University, 2008.

\bibitem{FFR}
Albert Fathi, Alessio Figalli, and Ludovic Rifford 
\newblock On the Hausdorff dimension of the Mather quotient.  
\newblock Comm. Pure Appl. Math.  62  (2009),  no. 4, 445--500.


%\bibitem{FGS}
%Albert Fathi, Alessandro Giuliani and Alfonso Sorrentino.
%\newblock Uniqueness of Invariant Lagrangian graphs in a homology or a cohomology class. 
%\newblock {\em Ann. Sc. Norm. Super. Pisa Cl. Sci.}, to appear.


%\bibitem{FathiSiconolfi}
%Albert Fathi and Antonio Siconolfi.
%\newblock Existence of {$C\sp 1$} critical subsolutions of the
  %{H}amilton-{J}acobi equation.
%\newblock {\em Invent. Math.}, 155(2):363--388, 2004.

%\bibitem{Katok}
%Anatole Katok and Boris Hasselblatt.
%\newblock Introduction to the modern theory of dynamical systems.
%\newblock {\em Encyclopedia of Mathematics and its Applications} Vol. 54, xviii+802,
%Cambridge University Press, 1995.


%\bibitem{Massart97a}
%Daniel Massart.
%\newblock Normes stables des surfaces.
%\newblock {\em C. R. Acad. Sci. Paris S\'er. I Math.},  324 (2): 221--224, 1997.



\bibitem{Massart97b}
Daniel Massart.
\newblock Stable norms of surfaces: local structure of the unit ball of rational
 directions.
\newblock {\em Geom. Funct. Anal.}, 7 (6): 996--1010, 1997.

\bibitem{Massart03}
Daniel Massart.
\newblock On Aubry sets and Mather's action functional.
\newblock {\em Israel J. Math.}, 134: 157--171, 2003.

\bibitem{Massart09}
Daniel Massart.
\newblock Aubry sets vs Mather sets in two degrees of freedom.
\newblock Preprint, 2010.


\bibitem{Mather90}
John~N. Mather.
\newblock Differentiability of the minimal average action as a function of the rotation number.
\newblock {\em Bol. Soc. Brasil. Mat. (N.S.)}, 21 (1):  59--70, 1990.
		

\bibitem{Mather91}
John~N. Mather.
\newblock Action minimizing invariant measures for positive definite
{L}agrangian systems.
\newblock {\em Math. Z.}, 207 (2): 169--207, 1991.

\bibitem{Mather93}
John~N. Mather.
\newblock Variational construction of connecting orbits.
\newblock {\em Ann. Inst. Fourier (Grenoble)}, 43 (5): 1349--1386, 1993.


\bibitem{Mather09}
John~N. Mather.
\newblock Order structure on action minimizing orbits.
\newblock {\em Symplectic topology and measure preserving dynamical systems}, 41--125, {\it Contemp. Math.} Vol. 512, Amer. Math. Soc., Providence, RI, 2010.




\bibitem{PPS}
Gabriel~P. Paternain, Leonid Polterovich and Karl~F. Siburg. 
\newblock Boundary rigidity for {L}agrangian submanifolds, non-removable
              intersections and {A}ubry-{M}ather theory.
\newblock {\em Mosc. Math. J.}, 3 (2): 593--619, 2003.



\bibitem{Sorrentino08}
Alfonso Sorrentino.
\newblock On the total disconnectedness of the quotient Aubry set.
\newblock {\em Ergodic Theory Dynam. Systems}, 28 (1): 267--290, 2008.

\bibitem{SorLecturenotes}
Alfonso Sorrentino.
\newblock Lecture notes on Mather's theory for Lagrangian systems.
\newblock Preprint 2010. (Available on ArXiv: 1011.0590).






\end{thebibliography}
\end{document}